# Dynamics and Control of Covid-19: Comments by Two Mathematicians

by Bernhelm Booß-Bavnbek (Roskilde, Denmark) and Klaus Krickeberg (Manglieu, France and Bielefeld, Germany)

Since the Covid-19 pandemic is not over it may appear to be premature to draw some conclusions. However, it may be, as well, just in time to recapitulate some lessons we as mathematicians should have learned and are urged to apply now. Thus we are asking: why are the dynamics and control of Covid-19 most interesting for mathematicians and why are mathematicians urgently needed for controlling the pandemic?

We shall first present our comments in a Bottom-up approach, i.e., following the events from their beginning as they evolved through time. They happened differently in different countries, and the main objective of this first part is to *compare* these evolutions in a few selected countries with each other.

Still, there are some general features, which we present separately as we are used to do in mathematics. They include the history of certain epidemics which have influenced the reactions of people in many countries, and some basic mathematical tools. In addition there is a common factor, which one of the present authors (KK) has defined on the 12th March 2020 in an e-mail to a German health office:

"The extension and evolution of Covid-19 in various countries and regions reflects the state of their health systems. This was for instance already very obvious in the case of Ebola."

It is in fact the *Public Health* component of the health system that plays a crucial role.

The second part of the article is not "country-oriented" but "problem-oriented". From a given problem we go "Top-down" to its solutions and their applications in concrete situations. We have organized this part by the mathematical methods that play a role in their solution. Here is an example where specially much mathematics is needed: to develop a vaccine and the strategy for applying it without loosing sight of basic ethical principles.

BOTTOM-UP

## 1. Prehistory

In the following the gentle reader may consult when necessary the book [KPP] for the basic concepts of epidemiology.

Demography as a mathematical subject area was already developed centuries ago well beyond its elementary beginnings. For a long time it remained the only mathematical tool in the study of the evolution of infectious diseases. Here is a famous early example. In China, India and Europe one tried to confer immunity against smallpox by infecting individuals slightly so they would contract a mild form of the disease and be immune afterwards. Some of them died by this procedure but in 1766 the Swiss mathematician Daniel Bernoulli showed by a demographical approach that the procedure would increase life expectancy if applied to everybody [DI1]. Nowadays evaluating the cost-effectiveness of a public health measure is being done widely; it is based on methods of mathematical economy.

The 19$^{th}$ Century saw the discovery of microorganisms as pathogens of many diseases and their study by mainly microbiological methods. The mathematical tools for following up an epidemic remained essentially demographical well into the 20$^{th}$ Century. A few physicians suggested that every epidemic ends because there are finally not enough people left to be infected, which is a naïve predecessor to the mathematical-epidemiologic concept of Herd Immunity (see Sect. 8). Nevertheless even the abundant literature on the influenza pandemic of 1918-19, wrongly called Spanish Flu, discusses only two possible ways for its ending: better clinical treatment and mutations of the pathogen.

Seen from a virological viewpoint the Spanish Flu was an extreme form of the so-called seasonal influenza. The virus which causes them can be one of a large variety, its genus being denoted by A, B, C or D, where some of them include several species. A is the most serious one; is has subtypes $A(H_xN_y)$, $x = 1,...,18$ and $y = 1,...,11$, where x and y represent proteins on the surface of the virus. The strategy for controlling the "normal" seasonal influenza epidemic is widely known even among laymen: identify the strain of the virus in autumn, develop a vaccine as fast as possible, and vaccinate people thought to be at risk. Nevertheless the number of infections and of deaths by a seasonal influenza can be as high as those by some of the pandemics to be described now.

The Spanish flu was due to $A(H_1N_1)$. Pictures from that time show people wearing masks that resembled those used now. In the years 1957-58 another "digression" from seasonal influenza occurred, called the Asian Flu and caused by $A(H_2N_2)$. It started in China and then became a pandemic, passing from neighbouring states through the UK and the USA. Estimations of the number of cases vary around 500 millions and of the number of deaths around 3 millions. Its beginnings looked much like those of the Spanish Flu but towards the end a vaccine became available, a predecessor to the ones being used now routinely against the seasonal flu.

The Hong Kong influenza of 1968-69, generated by the virus $A(H_3N_2)$, had similar characteristics and will not be described further.

Parallel to the entering the scene of these and other epidemics, and partly motivated by them, basically new mathematical tools of public health emerged in the first part of the 20[th] Century, preceded by a few studies in the late 19[th]. They were twofold. The first tool was called a "statistical-mathematical model". Its aim is the study of the influence of *factors*, also called determinants, on the health of people. Such factors may for instance be a lack of hygiene or a polluted environment. A factor can also be a preventive or curative treatment by an immunization or a drug, respectively; in that case the main objective of a study is to estimate the *efficacy* of the treatment. Sampling plans are statistical-mathematical models of a different but related kind. They form the basis of sample surveys, which are being done in profusion about Covid-19, too, and not always very illuminating.

The second tool is called "mathematical modelling of the evolution of an epidemic", or briefly "mathematical modelling". There are two kinds of it. First, one may aim at the epidemic curve, which is the cumulated number of cases up to a moment *t* as a function of *t*. In that case mathematical modelling serves to estimate or predict this curve under various assumptions on the infectivity of infected subjects. Early predecessors are presented in [FIN], see Figure 1; the question whether the infectivity remains constant or decreases played already a role. Refined versions are still being used, in particular for Covid-19 (Sect. 7).

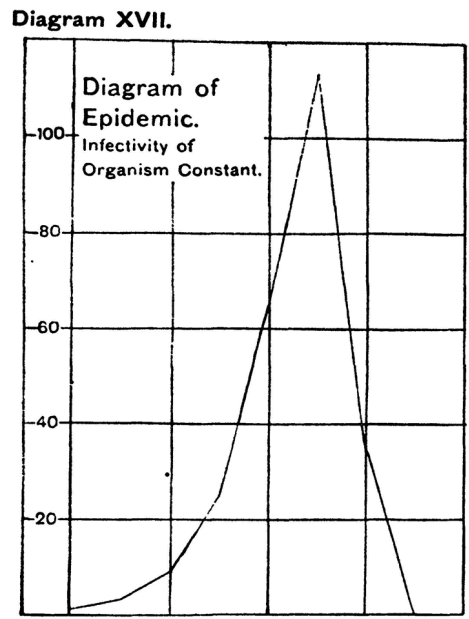

**Fig. 1** Early numerical simulation of an epidemic curve by J. Brownlee, 1907, discussed in [FIN]

Second, one may build so-called compartmental models (Sect. 8). The first one, for measles, was published in 1889 by P.D. En'ko; see [DI2]. Around the year 1900 compartmental models for malaria appeared. Then in the 1920s new models for the evolution of measles in closed populations were defined and intensively studied. They became very influential because they displayed already many basic features that reappeared later in mathematical models of epidemics in other and more complex settings.

Such tools found many applications. Dealing with large epidemics mathematically was no longer a matter of demography alone, although that continued to be the main tool for estimating number of cases and deaths. Statistical-mathematical models were employed to estimate the efficacy of antiviral drugs, for instance against HIV-infections, and the efficacy of various immunizations including those against forms of influenza. Mathematical modelling of epidemics was used in planning strategies to eradicate smallpox, poliomyelitis, measles and perhaps others. The first articles on modelling influenza epidemics appeared in the scientific literature. Planning a vaccination strategy involves both statistical-mathematical and mathematical models [HAl].

These roads to progress may have produced a general feeling of success in dealing with epidemics. Then in the period from 2002 to 2019 a few events occurred that evoked memories of previous pandemics and undermined such believes.

## 2. Unexpected Events 2002-2018

In November 2002 the first SARS (Severe Acute Respiratory Syndrome) epidemic broke out. It was a zoonosis generated by the virus SARS-CoV-1, a strain of the species SARS-CoV. It was first identified in China and never spread much beyond the surrounding countries and Canada. In July 2003 it was declared eradicated after having caused 8,096 cases and 774 deaths.

Similarly the Middle East Respiratory Syndrome (MERS), due to the corona virus MERS-CoV, lead to around 2,500 cases and 870 deaths between 2012 and 2020. It was essentially concentrated on Saudi Arabia and, to a minor degree, on South Korea, with most infections happening around the years 2014 to 2015. Being a zonoosis carried largely by camels, it is also called the Camel Flu.

Moreover a pandemic influenza invaded the world that resembled the Spanish Flu in several respects. Its pathogen was a new strain called $A(H_1N_1)09$ of the $H_1N_1$ influenza virus. Its origin is being debated; a likely hypothesis says that, being a zoonosis carried by pigs, it infected a human on a Mexican pig farm around January 2009. It was therefore called Swine Flu or Mexican Flu. It spread from North America to the whole world and was declared "extinguished" in August 2010. Estimations of the number of infections and deaths vary enormously but there were apparently more cases and fewer deaths than by the Spanish flu. Accusations against WHO were raised about its handling of vaccines against the Swine Flu.

Finally another zoonotic influenza appeared, popularly called Bird Flu and in scientific language Highly Pathogenic Avian Influenza (HPAI). The main pathogen was an $A(H_5N_1)$ influenza virus. It had been known long ago but reached a peak in the years 2013-2017. Whether there existed an airborne transmission from poultry to humans was a hotly debated question with obvious economic consequences. The Bird Flu spread widely over the whole world but the number of known human cases remained small, just over 70.

In addition to various forms of influenza and the epidemics generated by the corona virus SARS-CoV-1, SARS-CoV-2 or MERS-CoV, other epidemics occurred. It is instructive to compare them with those just mentioned, applying in addition mathematical yardsticks. We shall restrict ourselves to Ebola epidemics. Their most widespread outbreak was the Western African Ebola virus epidemic from 2013 to 2016, which caused 28,646 cases and 11,323 deaths. There is a fundamental difference in the evolution of a case of influenza or SARS-CoV-1 or SARS-CoV-2 on the one hand and of an Ebola case on the other, which leads to a basic difference in their mathematical modelling (Sects. 4 and 8). A carrier of an influenza or corona virus can transmit it to other persons well before the first symptoms appear, that is well before the end of the incubation period. A subject infected by Ebola will become infective only around the end of the incubation period. He (if it is a man) could then be immediately isolated together with his latest contacts in order to avoid further extension of the infection, provided that there is a health service nearby to do it. Therefore Ebola did not spread to countries that have a sufficiently dense primary health care network but it caused much suffering in countries that do not have it. The strategy of WHO to control the epidemic was wrong. It insisted on drugs and the search for a vaccine (which became available only in December 2019) but neglected primary health care. For the present purpose it would even have been most useful to rapidly train village health workers and "barefoot doctors" as it had been done decades ago.

## 3. Looking at some Countries

Only very few countries profited from the experiences of these premonitory 18 years to prepare much in advance for a possible, and probable, new outbreak of an epidemic. Some others took appropriate measures only at the first signs of Covid-19, and many started planning when the epidemic had almost reached its zenith.

We shall sketch some examples. For simplicity we shall always describe the result of the strategy of a country by indicating its cumulated numbers of confirmed cases and deaths around the 1st June 2020. Regarding the reliability of these data see Sects. 5 and 6.

We begin with those that had planned early.

**Taiwan:** Already in 2004, the year after the SARS-epidemic outbreak, the government established the National Health Command Center (NHCC), which was to prepare the country for a possible new epidemic. From 2017 on it was headed by the popular Minister of Health, Chen Shih-chung, who had studied dentistry at the Taipei Medical College. The Vice-President of Taiwan from 2016 to 2020, Chen Chien-jen, had been Minister of Health from 2003 to 2005 after having studied human genetics, public health, and epidemiology at the National Taiwan University and the Johns Hopkins University in the USA, followed by research. Thus decisions about the control of Covid-19 were taken by politicians competent in matters of health including Public Health.

Taiwan counts 23 million inhabitants and many of them travelled from and to China. From the 31st December 2019 on when WHO was notified of the epidemic in Wuhan all incoming flights from there were checked, followed by controls of passengers arriving from anywhere. An "Action Table" was produced in the period 20th January to 24th February 2020, which listed 124 measures to be taken. The public obtained daily revised clear information by all existing means. "Contact tracing", which means repeated follow-up of symptomatic persons, of confirmed cases and of all of their contacts, was rapidly established on the basis of the electronic health insurance card that everybody has. The virological PCR-tests used (Sect. 4) were already available and quarantines well organised. In late January rules about the wearing of masks were edited; a sufficient supply existed already.

As a result 442 confirmed cases had been found and 7 deaths recorded up to the 1st June.

**Vietnam:** The Vietnamese strategy resembles the Taiwanese one in almost all aspects, with the exception of contact tracing. A Steering Committee to deal with new epidemics existed within the Ministry of Health. It put into effect its plan right after the 23rd January when the first infected persons arrived at Vietnamese airports, among them a Vietnamese returning from the UK. All schools were closed on the 25th January, and since the 1st February everybody entering Vietnam must spend two weeks in quarantine.

Other measures were imposed or relieved in accordance with the evolution of the epidemic, for instance a limited confinement or the wearing of masks. The Ministry of Health issued regular precise and clear information for the entire population by all available means including smartphones. In addition there is a personalized information system by so-called "Survival Guides" given to everybody. Every survival guide defines three categories of persons: F0: a confirmed case; F1: suspected to be infected or having had contact with an infected person; F2: having had contact with a person in F1. Each person is expected to find the category to which it belongs. The survival guide then provides printed information about what she or he must do as a function of her or his category, for example to submit to a test. Only PCR-tests are being used.

In contrast to Taiwan contact tracing does not use electronic tools. It is being done by the population itself, aided by the survival guides, together with a large number of well-trained members of the health services, for example university lecturers.

At the end of 2019 Vietnam had 98,257,747 inhabitants. On the 1st of June there had been 328 confirmed cases and 0 deaths. These data are based on a strong demographic section of the "General Statistical Office" and on several Health Information Systems [KKR] and can hardly be contested.

The preceding sketch of control measures in Taiwan and Vietnam has shown us the three main components of their epidemiologic side: contact tracing; lock-down, that is physical, or social, distancing in the wide sense including quarantine and border controls; wearing of masks. We may call this the "surveillance-containment strategy". In addition there is the medical-clinical side, from primary health care such as general practitioners up to large hospitals. Its state is crucial to the number of deaths caused by the virus SARS-CoV-2.

In contrast to Taiwan and Vietnam it seems that all other countries of the world were unprepared at the end of December 2019. A few of them took fairly systematic and strict measures that covered the entire population as soon as the first cases had declared themselves. For a quick overview see Figure 2. This was for example true for **China** at the end of January 2020, for **Slovakia** and **Greece** on the 27th and 28th February, for **Austria** on the 10th March and for **Denmark** on the 12th March. An alternative Danish strategy, based on rigorous contact tracing and quarantine, but not implemented until now was argued for in [SIA].

Regarding the results the turbulent evolution in China is well known. In Denmark, with a population of 5.806 million, about 12,000 cases had been confirmed and 593 deaths recorded, and the corresponding figures for Austria are 8.86 million people, 16,979 cases and 672 deaths.

The confrontation of Slovakia, a country of around 5.5 million inhabitants, with Greece, which counts 10.72 million people, is particularly striking because it makes visible the role of their physicians and hospitals. In Slovakia there were 1,528 confirmed cases and 28 deaths. The corresponding data for Greece are 3,058 and 183. The relatively much higher number of fatalities in Greece, in spite of equally early reaction and almost the same number of cases per number of inhabitants, is no doubt due to the catastrophic state of its medical-clinical system caused mainly by the debt crisis from 2010 on.

Next we pass to a group of countries that reacted late and not systematically, applying the various measures in a haphazard way and only to part of the population. Here are some of them with their numbers of inhabitants in million, cumulated numbers of confirmed cases and numbers of fatalities:

**Table 1**. Countries that reacted late and not systematically (numbers pr. 1 June 2020)

| **Country** | Inhabitants (mio) | Confirmed cases | Fatalities |
|---|---|---|---|
| **Belgium** | 11.46 | 59,348 | 9,606 |
| **Spain** | 46.94 | 289,046 | 27,136 |
| **Italy** | 60.36 | 235,561 | 34,043 |
| **France** | 66.99 | 154,591 | 29,296 |
| **Germany** | 83.02 | 187,000 | 8,831 |

The relatively low number of deaths in Germany reflects mainly a sufficient medical-clinical system that could readily adapt itself to the epidemic. The opposite was true in France. There, about 100,000 hospital beds had been eliminated in the period between 1993 and 2018. An arbitrary strict "confinement" not determined by epidemiologic reasoning was imposed on the 17th March.

Finally there are countries that decided to do nothing, at least for a long while. Their motivation, or pretext, was above all a belief in herd immunity (Sects. 1 and 8) according to which the epidemic would stop by itself. This was the strategy of **Sweden**, a country counting 10.23 million persons, which resulted in 37,814 cases and 4,403 deaths. In the **United Kingdom** there were, among 66.65 million inhabitants, around 290,000 cases and 41,128 fatalities, and in the **USA** these data were 328.2 million, 2.04 million, 115,000 deaths.

This overview of strategies confirms that, as said in the introduction, the results depend indeed heavily on the state of Public Health. Note that nowadays in every language of the world the concept "Public Health" is designated by a literal translation or a slight modification of this expression. For instance in Danish it is "folkesundhed", that is, "Health of the People".

**TOP-DOWN**

In this second part we shall sketch the scientific and in particular mathematical principles involved in the study of successive stages of the pandemics. In short: Sect. 4: Discovery of the new virus, basic properties, testing for its presence in a person. 5 and 6: Data on the evolution of Covid-19 in a population. 7: Attempts at analysing mathematically and predicting such an evolution by representing it by an epidemic curve. 8: The analogous for a representation by a compartmental model. 9: Trying to stop the epidemic by a vaccine. 10: What to learn and what to do?

## 4. The New Virus SARS-CoV-2

After the often-depicted outbreak in late December 2019 of cases of pneumonia of unknown aetiology around Wuhan, in the course of January 2020 Chinese scientists identified a new virus as the pathogen. They followed the usual procedures, i.e., they determined the load of 26 common respiratory pathogens in the patients. They found none of them in abundance. They suspected SARS-CoV but could not find it either. Then they investigated all kinds of viral load that had a slight similarity (coincidence in a number of genomes) with SARS-CoV and detected a novel virus which displayed abundant virions in respiratory specimens from patients. Electron microscopy and mathematical pattern analysis [MUM, PEV] showed that it belongs to the same species as SARS-CoV-1 and MERS-CoV (Sect. 2); hence the name SARS-CoV-2.

Starting with this work in China a large number of publications about the peculiar properties of the pathogen and the ways it is acting have appeared. On the virological side its genetic sequence was determined. The new virus is believed to have zoonotic origins but human to human infection was rapidly established. The combination of SARS and Influenza features, that is intensive respiratory inhibition of patients and rapid transmission, make Covid-19, the disease caused by SARS-CoV-2, particularly dangerous. For further work see [AND].

In the clinical context, several periods in the evolution of a case were determined (see their definition in [KPP, Sect. 5.2]): The median *incubation period* is 5.2 days; the mean *latency period* is 4.6 days, i.e., in general the *infectious period* starts indeed before the *prodromal phase*. We have discussed the implications in Sect.

2 by comparison with Ebola. The mean length of the infectious period is 6 days for mild and asymptomatic cases; for severe and critical cases this period lasts on the average 22 days and ends only by recovery or death.

The manifold applications to the control of the pandemic of both their virological and their clinical characteristics will appear in Sects. 7, 8, 9 and 10. Their study is still active and may even reverse former results; this happened for example recently about so-called cross-immunities. However, in this article we shall only treat applications to the basic element of well-designed control strategies, namely testing for infections.

The first step of a test programme is to define the target population. Who will be tested? Subjects who had a contact with infected people? Or those who complain about symptoms? Or everybody coming from a region where cases exist? See the example of Vietnam in Sect. 3.

Next, what will be the objective? To discover the presence of the virus or that of some kind of antibodies? Depending on the objective there exist virological and serological tests. The usual virological test is called the PCR-(Polymerase Chain Reaction) Test. Dozens of serologic tests of varying quality have been and still are developed and even offered in some countries to the general public. Recall that the characterization of a test with a given target population and a given objective is a classical subject of clinical epidemiology [KPP, Sect. 19.2].

Coming back to the fundamental role of testing in control strategies we only remark that in poor countries or in rich countries with inattentive public health officials, the target population was often determined by the shortage of test kits and by the influence of institutions that required them for themselves.

## 5. Demography: Descriptive Epidemiology

This is classical *medical statistics*, which gives for a specific disease the number of cases and deaths together with the when and where and a few additional data such as sex, age and sometimes profession of the subjects.

In principle the methods for finding the number of confirmed cases and of fatalities by Covid-19 are the same as for any other disease. They fluctuate widely between countries. Both the diagnosis of a case of a disease and the description of the cause of a death may be relatively correct or most unreliable. In particular finding a correct diagnosis for somebody who complains about acute health problems depends very much on the local contact tracing methods and on the state of the clinical-medical system. An additional difficulty arises form the existence of asymptomatic forms of the disease, that is, subjects infected by SARS-CoV-2 who display no symptoms.

In Sect. 3 we have mentioned Vietnam, which uses its normal demographic and health information systems [KKR]. It includes in its statistics asymptomatic cases found by contact tracing. Other countries obtain their morbidity and mortality data from a "Health Reporting System". Such a system is partly based on sampling methods from various sources, for example hospitals and local health offices. In Germany the Robert Koch-Institute, a central institute mainly devoted to infectious diseases, reports on the results for Covid-19. In the USA the Johns Hopkins University plays a similar role. Still other countries use data from health insurance offices.

However, many countries have neither a health information system nor a health reporting system, or they do not use it for Covid-19. A host of alternative methods is being employed. For example France counts only hospitalized confirmed cases and only deaths which happen in a hospital or in a retirement home that is connected with a medical structure.

Summing up we may say that morbidity data and to a lesser degree mortality data for Covid-19 that one finds in various periodic publications are fairly unreliable, with very few exceptions. The sources are not always clearly indicated.

An important alternative idea is to compare the present situation with that in years past. Speaking again naïvely we assume that the present higher case frequencies and death tolls, and only these, are the result of Covid-19. Given the diagnostic difficulties mentioned above this idea is mainly applied to fatalities and hardly to nonlethal cases. Thus in the method of "excess mortality" we only measure how many more deaths by any cause happened this year than in the corresponding period in the past. For the UK we have for instance quoted in Sect. 3 the figure of 41,128 deaths up to the 1st June as supplied by the National Health Service. By contrast the National Statistical Office advanced about 62,000 deaths as excess mortality!

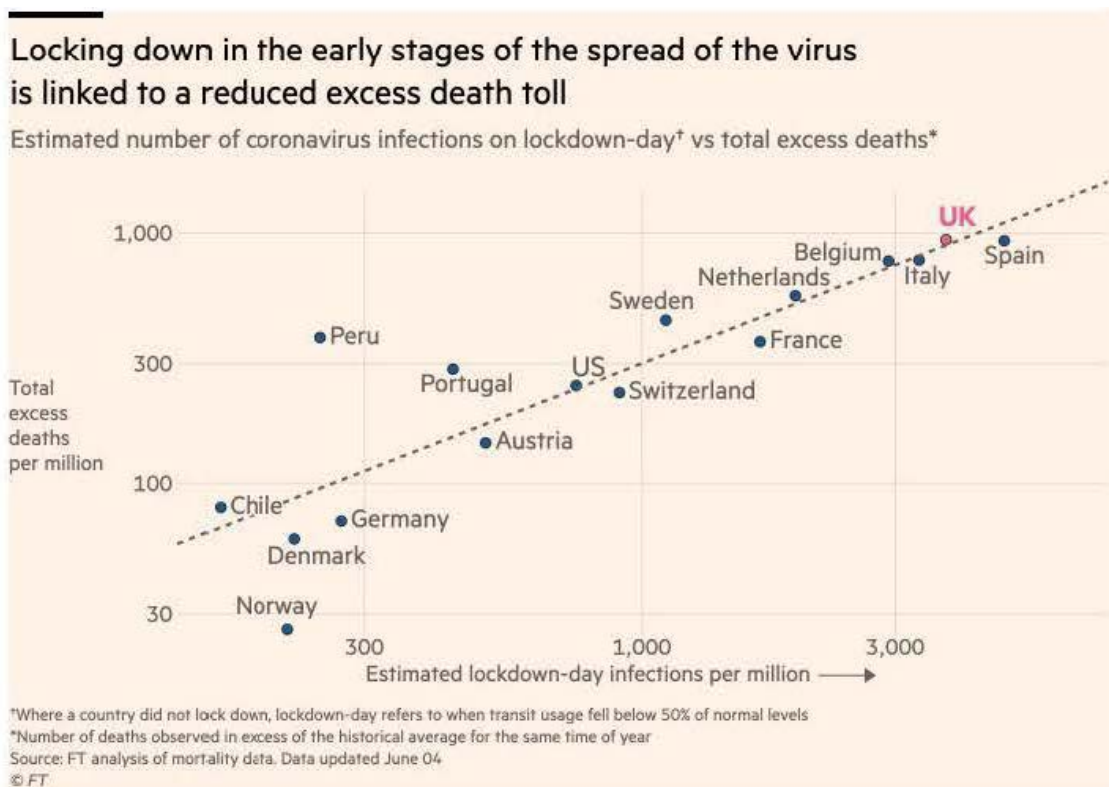

**Fig. 2** Estimated number of infections on lock-down day and excess mortality for 16 selected countries. Reproduced from [FIT], permission granted by *Financial Times Syndication*

Finally, here is an interesting idea based on the most classical form of a statistical-mathematical model. A graphic in the paper [FIT] (see Figure 2) shows for every one of 16 selected countries the point in the plane whose coordinates are, respectively, the estimated number of infections per million inhabitants on lock-down day and the excess mortality. A short glance convinces us that they are positively correlated. A simple regression analysis based on this graphic would allow us to estimate one of these values by the other one for any other country, too.

## 6. Advanced Demography

It goes in several directions beyond classical health statistics, all of them relevant to Covid-19, too. Firstly, *sample surveys* are conducted instead of using the data from the entire "target population". They have for example been used to study the influence of *social factors* on the evolution of various aspects of the disease. In particular the factor "to be an immigrant or to descend from them" was thoroughly investigated in some countries. Secondly, more types of data about cases and deaths are collected, for example about morbidity and mortality by age groups. Thirdly, data sets are not only being registered and perhaps published but also transformed and interpreted in various ways. Here, standardization is the best-known procedure. A fictitious example would be the number of fatalities by Covid-19 in Denmark if Denmark had the same age structure as Vietnam and in each age group it had the same Covid-19 mortality as in the same age group in Vietnam.

In Sect. 9 we shall meet statistical-mathematical models as a basic mathematical tool in developing a preventive treatment of Covid-19. With their help one studies in a clinical trial the influence of various factors on some outcome variable $E$ of interest. Here the idea of "controlling" for the influence of another factor, which might be a "confounder" in the study of the action of $E$, plays a role. It looks as if most demographers on the one hand, and most clinical epidemiologists on the other, ignore that the mathematical procedure of standardizing is the same as that of controlling for a confounder [KPP, Lesson 21]. A mathematician will not be astonished, though!

## 7. Modelling the Epidemic Curve

We have mentioned this classical concept in Sect. 1; see [KPP, Sect. 4.6]. Let $C$ be an epidemic, $V$ a geographical region, $t_0$ a moment of time which may be that of the first case of $C$ in $V$, and $f(t)$ for $t \geq t_0$ the number of observed and reported cases of $C$ that had declared themselves in $V$ before or at the instant $t$. Then $f$ is called the *epidemic curve* of $C$ in $V$. In particular it needs to be said whether unconfirmed cases are included or not. Measuring $f(t)$ as the time $t$ goes along is the task of the relevant demographic services (Sects. 5 and 6). This process is therefore subject to all the deficiencies listed there.

To get some knowledge about $f$ for various regions $V$ is of course one of the main concerns of the population of a country invaded by $C$. Such knowledge is equally vital for health authorities who attempt to control $C$. However, much more knowledge is desirable. What can we learn about the mechanism of $C$ by observing $f(t)$ ? This was already the subject of the papers described in [FIN]; see Sect. 1. In particular, is there a way to predict aspects of the future evolution of $f$, having observed the values $f(t)$ for a while?

Answers to these questions are generally given by *modelling f*, that is by making certain assumptions about its shape and by estimating certain parameters in it. A very large number of papers was published about this issue. Some of them use extrapolation methods known from mathematical economy. A recent survey on basic ideas and techniques can be found in [KRM] where a model is described in terms of an integro-differential equation.

We shall restrict ourselves to a discussion of an application, namely the so-called basic reproduction number $R_0$. It appears constantly in popular publications. To define it let us look at a subject $s$ that is infected at a time $t^* \geq t_0$. Let $\mu(s, t^*)$ be the number of all subjects infected by $s$ after $t^*$ in the form of secondary, tertiary etc. infections. Then $R_0$ is the average of $\mu(s, t^*)$ over all $s$. Thus it depends on $t^*$. It is precisely this dependence in which people are interested: a value less than 1 is looked upon as predictor of the extinction of $C$ after $t^*$. In the case $C$ = Covid-19, values as high as 5.7 had been estimated in the beginning, that is, for $t^*$ close to the time of the first outbreak of $C$. The article [SIA] presents an interesting factorisation of $R_0$ in order to compare different approaches to control the size of it.

## 8. Compartmental Models

We have sketched their historical origin in Sect. 1. We distinguished between two ways of mathematically modelling the evolution of an epidemic. Models of the first kind (Sect. 7) represent the temporal evolution of the number of subjects in a certain state, for instance the state "to be infected". By contrast, compartmental models also represent *changes* of this state at some moments in the form of transitions of a subject from one compartment to another one.

The SIR-model, which we designated in Sect. 1 as "intensively studied in the 1920s", is particularly simple and has served as a paragon for many others, in particular for those applied to Covid-19. It involves three compartments: *S* are the susceptible, not yet infected subjects, *I* the infected ones, and *R* consists of subjects removed by recovery with immunity or death. The transitions between compartments are described by differential equations for the numbers *S*(*t*), *I*(*t*) and *R*(*t*) of subjects in the compartments as a function of time *t*. They involve certain parameters such as transition probabilities from one compartment to another one. Under various assumptions the resulting system of differential equations for *S*, *I* and *R* can be solved explicitly or numerically.

A first important application is to estimate the basic reproduction number $R_0$ defined in Sect. 7. It can be expressed by the basic parameters.

Secondly, it turns out that the limit $S_\infty$ of *S*(*t*) for $t \to \infty$ is strictly positive, which means that a certain part of the population will never be infected. This led to the concept of *herd immunity*, which, however, gave rise to much confusion among people who thought they had something to say about the matter.

After the outbreak of Covid-19 many more involved compartmental models were defined and analysed. Their parameters represented among other features the underlying control strategy to be used. There was for instance the "do nothing" strategy and also the "mitigation" strategy, which consisted of the less stringent components of the "surveillance-containment strategy" defined in Sect. 3. In the much discussed paper [FER] Neil Ferguson and collaborators described the shape of the function *I*, that is the number of infected subjects, for the "do nothing" strategy. From the value 0 on it increases, reaches a maximum, decreases and finally reaches 0 at a certain moment $t_{happy}$. This had apparently motivated the countries UK, USA, Sweden and Brazil to adopt this strategy for too long, ignoring that Ferguson predicted (see Figure 3) about 500,000 deaths caused by the epidemic in the UK and 2.2 million in the USA before extinction at the moment $t_{happy}$.

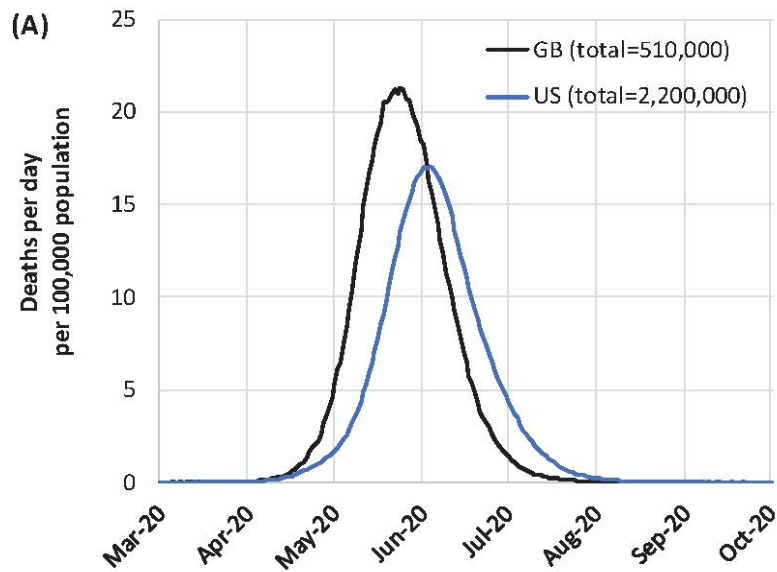

**Fig. 3** Expected deaths caused by the epidemic for the do-nothing strategy, reproduced from [FER] with permission of School of Public Health, Imperial College London

At present compartmental models play hardly any practical role, mainly because they contain too many unknown parameters. Some parameters such as infectivity are estimated with the help of a model of the epidemic curve, which seems to be a not very successful detour.

## 9. Preventive and Curative Treatments

It will hardly surprise that several pharmacological companies have started a run for developing curative and preventive treatments of various ailments, which SARS-CoV-2 may inflict on a person. Up to now no curative treatment was found. There are only the well-known methods to be used in the treatment of non-specific aspects of a case such as reducing pain, facilitate breathing or shorten the time to recovery by a antiviral drug. We shall therefore restrict ourselves to *preventive* treatments, that is, to immunizations.

The objective of an immunization by a vaccine against a Covid-19 connected health deficiency needs to be defined in the same way as for any other infectious disease. First the target population needs to be determined: whom do we intend to protect? Next, what are the health deficiencies we want to prevent? For how long is the preventive effect to last? This is a particularly important aspect of the vaccine but is usually suppressed when a new one is announced. For instance the measles vaccination remains lifelong active in most subjects. For Covid-19 the company which tries to develop the vaccine may be satisfied with a few months, hoping that SARS-CoV-2 will have disappeared after that. Finally the *efficacy* needs to be found, which represents the part of the target population actually protected. It may also be defined in epidemiologic terms by regarding as "exposed" all subjects that had not obtained the treatment. Then the efficacy is the "aetiological fraction among the exposed subjects".

Nowadays there is general agreement that the process of developing a vaccine against an infectious disease needs to run along a well-defined common line [KPP, Lesson 18, and HAL]. This ought to hold for Covid-19, too, and we shall therefore recall it here.

First, one or several substances are selected which, for some reasons whatsoever, usually virological ones, look like possible candidates for a vaccine. Each of them needs to be submitted to a "clinical trial" in order to explore its most important properties. Such a clinical trial consists of three "phases" I, II and III. Phase I deals with various mainly pharmacologic aspects such as side effects for various possible dosages.

Statistical-mathematical models are the essential tools of the phases II and III. Phase II aims at providing a first idea of the efficacy of the selected vaccine. Thus a relatively small target population is built artificially. Here two basic problems arise. The first is the definition of the outcome variable of interest. Often only the "immunogenicity" is being studied, which means the formation of antibodies, but not protection against the disease. It is a particularly complex and manifold problem in the case of Covid-19. Secondly, the target population needs to include among the vaccinated subjects a sufficient number of people who would attract the disease when not vaccinated. Since Covid-19 morbidity in the entire population of a country is small, such a group must be constructed by "challenge", that is, by infecting its members artificially. They are usually volunteers and their risk of dying is small except in the age groups where the lethality by the disease is high, that is, in the case of Covid-19, for old people. Faced with this ethical problem the USA used, for various previous infectious diseases, prison inmates whose terms were shortened as a reward. There was a time when Vietnam, while developing a certain vaccine, sent its samples for the Phase II trial to the USA to be tested in this way because Vietnamese ethical standards forbade all kinds of challenge.

There are usually several phase II trials in order to select the potential vaccine to be finally studied in a phase III trial. This is a field trial in the sense that a sample of subjects is drawn from the entire population of interest, for instance from among all inhabitants of a country within a certain age group. The outcome variable is not immunogenicity but protection against the disease in the sense of the desired efficacy. The size of the sample is determined beforehand by the precision of the intended estimate of the efficacy. As noted above the decision about the duration of the trial is a crucial element. If high efficacy during the first two weeks after vaccination is considered sufficient, the trial may be stopped after two weeks; this philosophy underlies the vaccinations against the seasonal influenza. If we are interested in its efficacy during the first ten years after vaccination, it must last ten years. This has, in addition to other problems, caused the long delay in developing an Ebola vaccine (end of Sect. 2). We hope that it will not be glossed over by those who are trying to sell a Covid-19 vaccine very soon.

## 10. Outlook

The pandemic has functioned like a magnifying glass. In some places, it showed a basically well-functioning society. In other places it revealed scandals and intolerable social inequalities. In particular it reflected the state of a country's public health system.

The present article aimed at describing the role of mathematics in the pandemic. As said above there are two parts to this "outlook". Let us take up the first one, namely: What can be learnt from the epidemic? In Sect. 1 we gave an overview of the main branches of mathematics that play a role. Then the Sects. 4 - 9 sketched the most frequent applications; their titles and their order correspond vaguely to the branches of mathematics concerned. Thus there were mathematical pattern analysis in laboratory work and statistical-mathematical models in judging the quality of tests; demographic methods in the collection of data; different ways to model the evolution of the pandemic mathematically; and clinical epidemiology in attempts to develop a vaccine.

In this way the article aimed at clarifying the potential role of mathematics in making decisions. On the one hand it turned out that in practise the role of epidemic curve or compartmental models is much more restricted than advertised in many publications. Decisions based on them may even have disastrous consequences, for instance those based on the mathematical concept of herd immunity. Thus blind trust in mathematical arguments is unjustified.

On the other hand denying the existence of a valid mathematical-scientific foundation for a control strategy is just as detrimental. It was done in Denmark with the "tracing and lock-down" strategy by a report of an "expert group" of health academics and officials, which reflected the interests of medical, industrial and governmental circles.

This comment leads us to the second part of our "outlook", namely: What to do in the future? The authors of the present article started it in early May by "Since the Covid-19 pandemic is not over ...". While we are finally finishing our work in the middle of July, it is still not over! It is even very active but has taken a largely different form. Hence it seems natural to analyse its present characteristics in the light of the facts we have described in the Sections 4 - 9 above and to ask ourselves: which lessons can we draw regarding the control strategies to be applied now?

Covid-19 does no longer surge from a single source. It reappears in small or large regions of many parts of the world, which may be of various forms and extensions: a single home for the elderly in France, two districts in Germany, a large city like Beijing, an entire province in Spain, or a whole country like New Zealand. We shall call them "nests" to distinguish them from "clusters", which denote certain discrete sets of people. A precise follow-up of the evolution of cases in these nests meets with the manifold difficulties explained in Sects. 5 and 6 and will not be repeated here.

A first natural question to ask is, then: why do "active" nests persist and reappear? Sect. 3 presented three components of successful control strategies: contact tracing; lock-down; masks. While contact tracing continues reluctantly, lock-down and wearing masks were widely abandoned, often as a result of governmental policies seeking popularity.

Next, what should be done? In the Sections 7, 8 and 9 we have explained, using in particular mathematical arguments, in how far the strategies of control treated there suffer from serious drawbacks. This leaves us with the combination of two measures: inside a nest a rigorous lock-down such as social distancing and preventing larger assemblies of people; at its borders: closing them or only allowing passage when combined with quarantine. For example New Zealand regarded as a single nest has taken such rigorous measures. As a result there are now no new cases, except two cases around the 14th July in "managed isolation facilities". Other nests will act similarly, we hope.

**Acknowledgements**:

Didier Dacunha-Castelle (Palaiseau, France),

Klaus Dietz (Tübingen, Germany)

**The authors**:

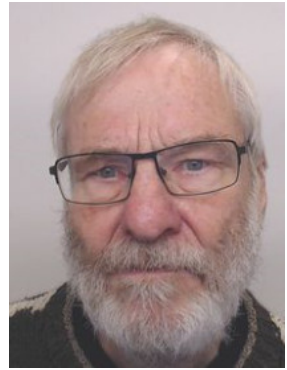

*Bernhelm Booß-Bavnbek* (booss@ruc.dk):

Born in 1941, studied mathematics from 1960 to 1965-at Bonn University. Research, teaching and practical work first in econometrics and operations research and then in geometric analysis and membrane processes of cell physiology. Affiliated to Roskilde University since 1977.

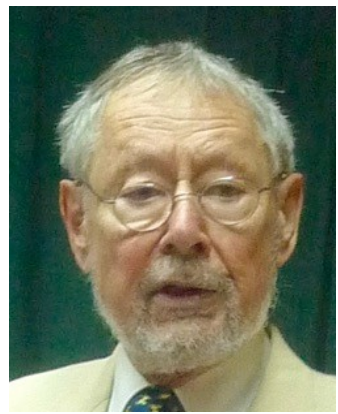

*Klaus Krickeberg* (krik@ideenwelt.de):

Born in 1929, studied mathematics from 1946 to 1951 at the Humboldt-University Berlin. Professor at several universities in Europe and outside; research, teaching and practical work first in mathematics and then in epidemiology and public health. Much of this was done in developing countries. Retired since 1998.